\documentclass[11pt]{amsart}

\usepackage{lmodern}
\usepackage{amsmath,amssymb,amsthm,mathtools,bm,hyperref}
\usepackage[margin=1.15in]{geometry}
\usepackage{microtype}
\usepackage{enumitem}
\usepackage{tikz}
\usepackage[capitalise]{cleveref}

\newtheorem{theorem}{Theorem}[section]

\newtheorem{lemma}[theorem]{Lemma}
\newtheorem{corollary}[theorem]{Corollary}
\newtheorem*{conjecture}{Conjecture}
\theoremstyle{definition}

\theoremstyle{remark}
\newtheorem*{remark*}{Remark}
\numberwithin{equation}{section}

\newcommand{\cc}[1]{\mathcal{ #1}}

\newcommand{\MI}{\mathrm{I}}
\newcommand{\KL}[2]{\mathrm{KL}\bigl(#1\,\big\|\,#2\bigr)}

\title{A proof of the weak Meyniel Conjecture}
\author{Joshua Erde}
\address{School of Mathematics, University of Birmingham, Birmingham, B15 2TT, UK}
\email{j.erde@bham.ac.uk}
\date{}

\begin{document}

\begin{abstract}
Cops and robbers is a pursuit game played on graphs by two players, one controlling a robber and one controlling a family of cops. The cop number $c(G)$ of a graph is defined to be the smallest number of cops which have a strategy to catch a robber. Probably the most well-known conjecture in the area is Meyniel's conjecture, that $c(G) = O( \sqrt{n})$ for every connected $n$-vertex graph $G$. Despite much attention, after 40 years the best known bounds are of the form $c(G) \leq n^{1-o(1)}$ and it is an open problem, known as the weak Meyniel conjecture, to give a bound of the form $n^{1-\varepsilon}$ for some $\varepsilon >0$.

In this paper we give a proof of the weak Meyniel conjecture, showing that $c(G) = O\left( n^{\frac{3}{4}}\ln n\right)$. 
Our proof is probabilistic in nature, relying on an information theoretical argument to show that for every fixed distribution on robber walks of length $\textrm{diam}(G)\cdot \ln n$, there is a deterministic adaptive strategy for one cop that intercepts a walk sampled from this distribution with probability $\Omega\left(n^{-\frac{1}{2}}\right)$.
\end{abstract}
\maketitle

\section{Introduction}
The game of \emph{cops and robbers} is a pursuit game played on a connected graph $G=(V,E)$ by two players with full information, one player controlling a set of $k \geq 1$ \emph{cops} and the other controlling a \emph{robber}. To start the game, the first player (who we will just refer to as the cops) chooses a starting vertices for each of the $k$ cops, and then the second player (who we will refer to as the robber) chooses an starting vertex for the robber. In each subsequent round, both players make a move, starting with the cops. Each cop may move to a vertex adjacent to his current location, or stay still, and the same goes for a subsequent move of the robber. The cops win if at some point there is a cop on the same vertex as the robber, otherwise the robber wins. This game was introduced independently by Nowakowski and Winkler \cite{NW83} and Quilliot \cite{Q78}.

The \emph{cop number} $c(G)$ of a graph $G$ to be the smallest number of cops $k$ such that the first player has a winning strategy in this game. Perhaps the most well-known open problem in this area is \emph{Meyniel's conjecture} (see \cite{F87}).

\begin{conjecture}[Meyniel's conjecture]
If $G$ is a connected $n$-vertex graph, then $c(G) = O(\sqrt{n})$.
\end{conjecture}

There has been much interest in this problem, see for example the survey of Baird and Bonato \cite{BB12} on the conjecture, where it is called one of the deepest open problems on the cop number of a graph. However, despite this, relatively little progress has been made on the trivial upper bound $c(G) \leq n$. Following results of Frankl \cite{F87} and Chiniforooshan \cite{C08}, currently the best known bound is
\[
c(G) = O\left( \frac{n}{2^{(1-o(1))\sqrt{\log{n}}}} \right) = O\left(n^{1-o(1)}\right).
\]
proved independently by Lu and Peng \cite{LP12} and Scott and Sudakov \cite{SS11}, and generalised to a more general class of pursuit games by Frieze, Krivelevich and Loh \cite{FKL12}.

If true, Meyniel's conjecture is best possible up to the leading constant, since it follows from a result of Aigner and Fromme \cite{AF84} that the incidence graphs of projective planes satisfy $c(G) = \Omega(\sqrt{n})$ (see \cite[Section 3]{BB12} for details).

Meyniel's conjecture has been shown to hold for a number of restricted graph classes, such as graphs of diameter two and bipartite graphs of diameter three \cite{LP12}, binomial random graphs \cite{PW16}, random regular graphs \cite{PW19} and most recently abelian Cayley graphs \cite{B20,BHT21}.

Due to the perceived difficulty of the conjecture, and the gap between the upper and lower bounds, the so called \emph{weak Meyniel conjecture} (or soft Meyniel conjecture) has gained prominence, which asks for an upper bound of the form $c(G) = O\left(n^{1- \varepsilon}\right)$ for some $\varepsilon >0$.

The main purpose of this note is to prove the weak Meyniel conjecture, in particular we show that $c(G) = O\left(n^{1-\varepsilon}\right)$ for any fixed $0 \leq \varepsilon < \frac{1}{4}$.

\begin{theorem}\label{t:weakmeyniel}
Let $G$ be a connected $n$-vertex graph with maximum degree $\Delta$, then $c(G) = O\left( n^{\frac{3}{4}} \sqrt{\ln n \ln (\Delta +1)}\right)$.
\end{theorem}

We will prove \Cref{t:weakmeyniel} using a mixed bound on the cop number in terms of the diameter and maximum degree of $G$.
\begin{theorem}\label{t:copdiambound}
Let $G$ be a connected $n$-vertex graph with diameter $D$ and maximum degree $\Delta$, then
\[
c(G) = O\left( D \sqrt{n} \ln n  \ln (\Delta +1) \right).
\]
\end{theorem}
Indeed, a standard argument dating back to Frankl \cite{F87} allows us to use $\frac{n}{d+1}$ cops to restrict the robber to a subgraph of diameter at most $d$ to which we can apply \Cref{t:copdiambound}. Optimising the choice of $d$ leads to \Cref{t:weakmeyniel}.

Another interesting corollary of \Cref{t:copdiambound} is that it implies that Meyniel's conjecture holds up to a polylogarithmic factor for graphs with small diameter.
\begin{corollary}\label{c:lowdiameter}
Let $G$ be a connected $n$-vertex graph whose diameter is polylogarithmic in $n$, then $c(G) = \tilde{O}(\sqrt{n})$, where the tilde hides a polylogarithmic factor.
\end{corollary}
This in particular improves results of Lu and Peng \cite{LP12} and Hosseini, Knox and Mohar \cite{HFM20} on the cop number of graphs of fixed diameter and of Bradshaw, Hosseini, Mohar and Stacho \cite{BHMS23} on the cop number of bounded degree expanders.

\subsection{Methods}
Let us sketch the proof of the key result, \Cref{t:copdiambound}.
We start by fixing $T=\Theta(D \log n)$ and only consider the first $T$ rounds of the game. Let $\cc{W}_T$ be the set of all possible \emph{trajectories} for the robber during this time period - i.e., each $\bm{r} = (r_0,r_1,\ldots, r_T) \in \cc{W}_T$ is a (lazy)-walk of length $T$ in the graph $G$.

The first key idea will be to think about playing against a random robber, whose trajectory $R$ is drawn from some known distribution on $\mathcal{W}_T$, which might be supported on some strict subset of $\mathcal{W}_T$. We will give a randomised adaptive strategy for a single cop which will have the property that there is some time $s \in [T]$ such that the probability that the cop catches the sampled walk $R$ at time $s$ is at least $p=p(n)=\Omega\left( n^{-\frac{1}{2}}\right) $. In particular, if $R$ is uniformly distributed on some subset $\cc{W}' \subseteq \cc{W}_T$, then we can conclude that there is a deterministic adaptive cop strategy which intercepts at least a $p$ proportion of the possible trajectories in $\cc{W}'$, in that any robber who followed that trajectory would have been caught. 

In this manner we can sequentially construct a series of cop strategies, each of which excludes a $p$ proportion of the remaining trajectories until none remain, which would require $\left\lceil \frac{1 + \log |\mathcal{W}_T|}{p} \right \rceil$ cops, which leads to the bound in \Cref{t:copdiambound}.

So, the crux of the matter is how to play against the random $R$-robber. The key new idea will be to try to intercept the robber in both time and space. The cop will maintain a (random) prediction $Y_s$ for the robber's location at a target time $b$ and arrange to be at predicted vertex at an earlier time $a \leq b$. In each subsequent round, the cop will move his target time one round earlier and update his prediction, taking into account the robber's subsequent moves, moving to this new prediction. In this manner, as the robber moves forward in time, the cop moves backwards along his `likely path' to intercept him. Let us write $\mu_s$ for the distribution of $Y_s$ and write $\rho_s$ for the \emph{forecast} of where the robber will be - the actual distribution of the robber's location at the target time given his observed trajectory so far.

We will be able to bound the difference between $\rho_s$ and $\mu_s$ in terms of information contained in the robber's moves about his later location during this interception phase (see \Cref{l:capture prob}, in particular \Cref{e:distdiff}). At the very end of the process, when the cop is predicting the actual time at which the robber is temporally located, the only difference between the two distributions will be the cumulative effect of this revealed information. In particular, since at this point the forecast is just a point mass at the robber's actual location, we will be able to use a bound on this cumulative information to bound the probability that the cop intercepts the robber at this target time (see \Cref{l:abscheme}). An averaging argument will show that for some choice of $a$ and $b$, relatively little information about the robber's future location is revealed in this time period, and so the probability of interception is quite large (see \Cref{l:capture}).

\begin{remark*}
Essentially the same bound on $c(G)$ as in \Cref{t:weakmeyniel} was also obtained independently by Geneson \cite{G26} using similar methods.
\end{remark*}

\subsubsection*{Declaration of AI usage}
The initial proof strategy was developed using ChatGPT 6.1 Sol, and it and Claude Opus 5.5 were used to assist in the drafting of the paper. The author takes full responsibility for the mathematical correctness and exposition of the manuscript.

\section{Preliminaries}
Let $G$ be a graph. Given $x \in V(G)$ we write $N(x)$ for the \emph{neighbourhood} of $x$ and $N[x]:= \{x\} \cup N(x)$ for the \emph{closed neighbourhood} of $x$. Given $T \in \mathbb{N}$, we define $\cc{W}_T =\mathcal{W}_T(G)$ to be the set of all lazy-walks of length $T$ in $G$, that is, all sequences $\bm{r}=(r_0,\ldots,r_T)$ where $r_i \in N[r_{i-1}]$ for each $i\in [T]$. In what follows, we will only talk about lazy-walks, which we will refer to as walks for convenience.

Given a random vector $R=(R_0, R_1, \ldots, R_T)$ and $s \in \{-1,0,1,\ldots,T\}$ we denote by $R_{\leq s} := (R_0,R_1,\ldots, R_s)$, with the convention that $R_{\leq -1}$ is empty. We write $\mathbb{P}_s$ for the conditional probability measure given $R_{\leq s}$. Given a random variable $X$ which is jointly distributed with $R$ we write $\mathbb{P}_s(X \in \cdot)$ for the conditional law of $X$ given $R_{\leq s}$. Note that $\mathbb{P}_s$ and $\mathbb{P}_s(X \in \cdot)$ are random variables and not a fixed probability measure or law.

Let $X,Y$ and $Z$ be jointly distributed discrete random variables with probability mass functions $p_X, p_Y$ and $p_Z$, respectively. We write $H(X) := - \sum_{x} p_X(x) \log_2 p_X(x)$, where terms with $p_X(x) =0$ contribute zero to the sum. This is the \emph{(Shannon) entropy} of $X$. The conditional entropy of $X$ given $Y$ is $H(X\mid Y) = \sum_{y} p_Y(y)H(X \mid Y=y)$ and the \emph{conditional mutual information} of $X$ and $Y$ given $Z$ is defined as $\MI(X;Y\mid Z) = H(X\mid Z) + H(Y\mid Z) - H(X,Y\mid Z)$, where we write $\MI(X;Y) = \MI(X;Y\mid \emptyset)$. We will write $\MI_s(X;Y)$ for the mutual information under $\mathbb{P}_s$. In the paper we will consider quantities of the form $\MI(R_a ; R_b \mid R_{\leq s})$ and $\MI_s(R_a ; R_b)$; let us comment briefly on the difference between them, which is subtle. The former is a real number, whereas the latter is a random variable, and they are related by the equality 
\[
\MI(R_a ; R_b \mid R_{\leq s}) = \mathbb{E}(\MI_s(R_a ; R_b)).
\]
We will require the following elementary facts about mutual information, all of which can be found, for example, in Cover and Thomas \cite{CT91}. 
Firstly, if $Z$ determines either $X$ or $Y$, then $\MI(X;Y\mid Z)=0$ and if $Z$ determines $X$, then $\MI(X;Z) = H(X)$. Secondly, the fact that a uniform distribution maximises entropy, and so in particular if $X$ is supported on a set of size $n$, then $H(X) \leq \log_2 n$. Finally, a simple consequence of what is known as the \emph{chain rule for mutual information}, together with the two facts above, lead to the following:
\begin{equation}\label{e:chainrule}
\sum_{s=0}^{t} \MI(R_s ; R_t \mid R_{\leq s-1}) = \MI(R_{\leq t} ; R_t) = H(R_t) \leq \log_2 n.
\end{equation}

Finally, given two distributions $\rho$ and $\mu$ defined on the same finite set $\cc{X}$ the \emph{Kullback-Leibler divergence} (KL-divergence) is defined as
\[
\KL{\rho}{\mu} = \sum_{x\in \cc{X} } \rho (x) \log_2 \left( \frac{\rho(x)}{\mu(x)} \right) \in [0,\infty],
\]
where terms with $\rho(x) =0$ contribute $0$ to the sum and the divergence is infinite if $\rho(x) > 0 = \mu(x)$ for some $x$.

The following identity is a consequence of the chain rule for KL-divergence; see Cover and Thomas
\cite[Theorem~2.5.3]{CT91}.

\begin{lemma}\label{l:mixture}
Let $\rho_1,\ldots,\rho_k$ and $\mu$ be probability distributions on a finite set, and let $\lambda_1,\ldots,\lambda_k \geq 0$ satisfy $\sum_{i=1}^k\lambda_i=1$. Writing $\rho=\sum_{i=1}^k\lambda_i\rho_i$, we have
\[
\sum_{i=1}^k \lambda_i \KL{\rho_i}{\mu} = \sum_{i=1}^k\lambda_i \KL{\rho_i}{\rho} + \KL{\rho}{\mu},
\]
where terms with $\lambda_i=0$ contribute $0$ to either side.
\end{lemma}
Note that, \Cref{l:mixture} holds with both sides taking values in $[0,\infty]$. When we apply the above we will use the following fact relating mutual information and KL-divergence, which can be easily verified from their definitions.

\begin{lemma}\label{l:MIKL}
Let $X$ and $Z$ be jointly distributed discrete random variables. Let $\rho$ be the law of $X$ and let $\rho_z$ be the conditional law of $X$ given $Z=z$. Then
\[
\MI(X;Z) = \sum_z p_Z(z) \KL{\rho_z}{\rho}.
\]
\end{lemma}

We also use the data-processing inequality for KL-divergence;
see Csisz{\'a}r and K{\"o}rner, \cite[Lemma 3.11]{CK15}

\begin{lemma}\label{l:dataprocess}
Let $\rho$ and $\mu$ be probability distributions
on a finite set $V$, and let $K$ be a stochastic
matrix on $V$. Then
\[
\operatorname{KL}(\rho K\Vert\mu K)
\leq \operatorname{KL}(\rho\Vert\mu).
\]
\end{lemma}

\section{Catching a random robber}
Suppose that we have a random vector $R=(R_0,R_1,\ldots,R_T)$ distributed on the set of walks $\cc{W}_T$ of length $T$ in a graph $G$ which represent the possible trajectories of a random robber. A \emph{random $R$-robber} is a robber whose strategy is to start at $R_0$ and follow the walk given by $R$, so that the robber's position after their move in the $s$th round is $R_s$.

We write $C_0$ for the cop's starting vertex and $C_s$ for their position after their move in the $s$th round, so that $C_s$ is chosen knowing $R_{\leq s-1}$. $C_s$ will be a function of both the random variable $R_{\leq s -1}$ as well as some additional randomness which is independent of $R$. Note that the robber is caught in round $0$ if $C_0 = R_0$ and caught in round $s \geq 1$ if $C_s \in \{ R_{s-1},R_s\}$, in which case we say that the cop \emph{intercepts} the trajectory $R$.

We will define a collection of randomised strategies for individual cops, each of which is determined by two integers $-1 \leq s_0 \leq s_1 \leq T$ and a sequence $\tau=(\tau_s \colon  s \in [s_0,s_1])$ of \emph{target times} $\tau_s \in \{0,1,\ldots,T\}$, where $|\tau_{s+1} - \tau_s| \leq 1$ for all $s \in [s_0,s_1-1]$. Let us refer to the triple $(s_0,s_1,\tau)$ as a \emph{scheme}. 

We will define a sequence of random variables $Y=(Y_s \colon s \in [s_0,s_1])$, jointly distributed with $R$, which represents the cop's \emph{predictions} for the location of the robber at the target times $\tau_s$. Later, we will give a strategy for a cop which guarantees that $C_s=Y_s$ for some terminal segment of the $s \in [s_0,s_1]$.

For all $s > s_0$, we will think of the cop as revealing the value of the random variable $Y_s$, which will depend only on $R_{\leq s-1}$ together with other independent random choices, at the start of the round $s$ before the cop makes his move. For technical reasons, the initial prediction $Y_{s_0}$ is revealed after the robber's move in round $s_0$, or before the initial placement if $s_0=-1$. Any further randomness used in defining $Y_s$ is drawn independently from $R$ and from $Y_{\leq s - 1}$. $Y$ is then determined by the following process
\begin{enumerate}[label=(Y\arabic*)]
    \item\label{i:initialY} $Y_{s_0}$ is sampled from $\mathbb{P}_{s_0} \left( R_{\tau_{s_0}} \in \cdot \right)$;
    \item\label{i:updateY} For $s \in [s_0,s_1-1]$, $Y_{s+1}$ is sampled from $K_s(Y_s,\cdot)$, where
    \[
    K_s(x,y) = \mathbb{P}_s\left( R_{\tau_{s+1}} =y  \mid R_{\tau_s} = x \right) \text{ if } \mathbb{P}_s\left( R_{\tau_s} = x \right) > 0,
    \]
    and $K_s(x,y) = \mathbf{1}_{y=x}$ otherwise.
\end{enumerate}
In other words, at the start $Y_{s_0}$ has the same distribution as the robber's location at the target time $\tau_{s_0}$ conditioned on the observed history $R_{\leq s_0}$. Then, in each update step, we sample $Y_{s+1}$ from the conditional distribution of the robber's location at the new target time $\tau_{s+1}$ conditioned on the observed history $R_{\leq s}$ in the previous round, conditioning further on the event that the previous value $Y_s$ was the correct prediction for the robber's location at the previous target time $R_{\tau_s}$.

Let us note two immediate properties of the transition matrix $K_s$:
\begin{enumerate}[label=(K\arabic*)]
\item\label{i:cstpred} If $\tau_{s+1} = \tau_s$, then $K_s(x,y) =  \mathbf{1}_{y=x}$ and so $Y_{s+1} = Y_s$;
\item\label{i:feasible} Since $R$ is a walk in $G$ and $|\tau_{s+1} - \tau_s| \leq 1$, it follows that $K_s(x, \cdot)$ is supported on $N[x]$.
\end{enumerate}

Let us define the \emph{information matrix} of $R$ as the matrix $\Gamma$ where
\[
\gamma_{s,t} = \MI( R_s ; R_t \mid R_{\leq s-1} ) \text{ for $0 \leq s,t \leq T$ }
\]
with the convention that $R_{-1}$ represents the empty conditioning. Let us note some elementary properties of $\Gamma$:
\begin{enumerate}[label=(A\arabic*)]
\item\label{i:lowerdiagonal} If $s > t$, then $\gamma_{s,t} = 0$;
\item\label{i:columnsum} By \eqref{e:chainrule}, for any $0 \leq t \leq T$, 
\[
\sum_{s=0}^T \gamma_{s,t} = \sum_{s=0}^T \MI( R_s ; R_t \mid R_{\leq s-1} )  \leq \log_2 |V|.
\]
\end{enumerate}

Let us start by bounding from below the probability that the guess $Y_{s_1}$ for the location $R_{\tau_{s_1}}$ of the robber at the target time $\tau_{s_1}$ is correct, in terms of the information matrix. To this end, we define 
\begin{equation}\label{e:weightdef}
W(\tau) =\sum_{s=s_0+1}^{s_1} \gamma_{s,\tau_s}.
\end{equation}

\begin{lemma}\label{l:capture prob}
Whenever $\tau$ is such that $\tau_{s_1} \leq s_1$, it follows that
\[
\mathbb{P}(Y_{s_1} = R_{\tau_{s_1}}) \geq 2^{-W(\tau)}.
\]
\end{lemma}
\begin{proof}
For each $s \in [s_0,s_1]$, let us define two distributions. The first $\rho_s = \mathbb{P}_s\left( R_{\tau_s} \in \cdot \right)$ is the distribution of the robber's location at the current target time and the second $\mu_s = \mathbb{P}_s(Y_s \in \cdot)$ is the distribution of the cop's prediction (of the robber's location at the target time) at time $s$.  

We will keep track of the expected KL-divergence $\Phi_s := \mathbb{E}(\KL{\rho_s}{\mu_s})$ between these two distributions. Let us first note a few simple facts about $\Phi_s$ and how the distributions $\rho_s$ and $\mu_s$ evolve.

Firstly, since $\mu_{s_0} = \rho_{s_0}$ by \ref{i:initialY}, we have that $\Phi_{s_0}=0$. Secondly, by \ref{i:updateY}, $\mu_{s+1} = \mu_s K_s$ and in particular $\mu_{s+1}$ is a function of $R_{\leq s}$. Indeed, $Y_{s+1}$ is determined by $R_{\leq s}$ together with randomness independent of $R$, and so revealing $R_{s+1}$ does not change its conditional law. Finally, by the law of conditional total probability
\begin{align*}
(\rho_s K_s)(y)  &= \sum_{x } \rho_s(x)K_s(x,y) \\
&= \sum_{x }\mathbb{P}_s\left( R_{\tau_s} =x \right)\mathbb{P}_s\left( R_{\tau_{s+1}} =y  \mid R_{\tau_s} = x \right) \\
&= \mathbb{P}_s \left( R_{\tau_{s+1}} = y \right),
\end{align*}
and so $\rho_s K_s = \mathbb{P}_s \left( R_{\tau_{s+1}} \in \cdot \right)$. 

Note that $\rho_s K_s$ is not quite $\rho_{s+1} =  \mathbb{P}_{s+1} \left( R_{\tau_{s+1}} \in \cdot \right)$. However, we can control how the divergence changes by substituting one for the other using \Cref{l:mixture}. We apply the lemma working with the probability measure $\mathbb{P}_s$ and taking $\lambda_x = \mathbb{P}_s(R_{s+1}=x)$, $\rho_x = \mathbb{P}_s(R_{\tau_{s+1}} \in \cdot \mid R_{s+1}=x)$, $\rho = \rho_s K_s$ and $\mu = \mu_{s+1} =\mu_s K_s$. Note that
\begin{align*}
\sum_{x\in V(G)} \lambda_x \KL{\rho_x}{\mu} &= \sum_{x \in V(G)}\mathbb{P}_s(R_{s+1}=x) \KL{\mathbb{P}_s(R_{\tau_{s+1}} \in \cdot \mid R_{s+1}=x)}{\mu_{s+1}}\\
&= \mathbb{E}( \KL{\rho_{s+1}}{\mu_{s+1}} \mid R_{\leq s}).
\end{align*}
Furthermore, by \Cref{l:MIKL}
\[
 \sum_{x \in V(G)}\lambda_x \KL{\rho_x}{\rho} = \sum_{x \in V(G)}\mathbb{P}_s(R_{s+1}=x)  \KL{\mathbb{P}_s(R_{\tau_{s+1}} \in \cdot \mid R_{s+1}=x)}{\rho_s K_s} = \MI_s(R_{\tau_{s+1}} ; R_{s+1} ).
\]
Finally, by \Cref{l:dataprocess}
\[
\KL{\rho}{\mu}  = \KL{\rho_s K_s}{\mu_s K_s} \leq \KL{\rho_s}{\mu_s}.
\]

Hence, we conclude by \Cref{l:mixture} that
\begin{align}\label{e:distdiff}
\mathbb{E} \left( \KL{\rho_{s+1}}{\mu_{s+1}} \mid R_{\leq s} \right)\leq \MI_s(R_{\tau_{s+1}} ; R_{s+1} ) + \KL{\rho_{s}}{\mu_{s}}.
\end{align}
Taking expectations of both sides, and using that 
\[
\mathbb{E}(\MI_s(R_{\tau_{s+1}} ; R_{s+1} ) = \MI(R_{s+1}; R_{\tau_{s+1}} \mid R_{\leq s}) = \gamma_{s+1,\tau_{s+1} },
\]
gives
\begin{equation}
\Phi_{s+1} \leq \gamma_{s+1,\tau_{s+1}} + \Phi_s,
\end{equation}
and so $\Phi_{s_1} \leq \Phi_{s_0} + \sum_{s= s_0}^{s_1-1} \gamma_{s+1,\tau_{s+1}} = W(\tau)$.

Finally, since $\tau_{s_1} \leq s_1$, the forecast $\rho_{s_1}$ of the robber's location at time $\tau_{s_1}$ given the observed history $R_{\leq s_1}$ is just a point mass at $R_{\tau_{s_1}}$. Hence, 
\[
\mathbb{P}_{s_1}\left(Y_{s_1} = R_{\tau_{s_1}}\right) = \mu_{s_1}\left(R_{\tau_{s_1}}\right) = 2^{-\KL{\rho_{s_1}}{\mu_{s_1}}}.
\]
and applying Jensen's inequality to the convex function $2^{-t}$ and taking expectations over $R_{\leq s_1}$ gives
\[
\mathbb{P}\left(Y_{s_1} = R_{\tau_{s_1}}\right) \geq 2^{- \Phi_{s_1}} \geq 2^{-W(\tau)},
\]
as claimed.
\end{proof}

Let us now define our family of cop strategies, which will be indexed by pairs of integers $(a,b)$ with $0 \leq a\leq b \leq T$. Roughly, the cop's strategy will be to intercept the predicted location $R_b$ at time $a$ and then follow this prediction backwards.

To this end we define the \emph{$(a,b)$-scheme} to be the triple $(s_0,s_1,\tau)$ where
\begin{itemize}
\item $s_0 := \max \{a-D,-1\}$;
\item $s_1 =  \lceil \frac{a+b}{2} \rceil$;
\item $\tau_s = \min \{ b, a +b -s\}$ for $s \in [s_0,s_1]$.
\end{itemize}
In other words, the target time is $b$ until round $a$, at which point it decreases by one per round. 

If $a \geq D$, then the cop will spend the $D$ rounds $a-D+1, a-D+2, \ldots, a$ travelling to $Y_{s_0}$, the predicted location (under $\mathbb{P}_{s_0}$) of the robber at time $b$, and from there track him backwards in time. Since $a \geq D$, the cop will always have time to reach $Y_{s_0}$ in this period. At first glance, it might seem that choosing $a < D$ poses a problem, as the cop might not have time to locate himself at the prediction. However, this is in fact simpler for us, since by \ref{i:initialY} the initial prediction $Y_{-1}$ is revealed before the cop chooses his initial location, and the cop can simply choose to start with $C_0 = Y_{-1}$ and stay there until time $a$.

Let us write $W(a,b) = W(\tau)$ and note that $W(a,b)$ splits naturally into two parts
\begin{equation}\label{e:locationandtrackingcost}
W(a,b) = \sum_{s=s_0+1}^a \gamma_{s,b} + \sum_{s=a+1}^{s_1} \gamma_{s,a+b-s},
\end{equation}
where the first sum represents the \emph{placement cost}, the cost of physically getting the cop to the predicted vertex and the second sum represents the \emph{interception cost}, the cost of tracking the robber backwards through time, where these costs contribute to the eventual KL-divergence between the distribution of the cop and the robber at time $s_1$.

\begin{lemma}\label{l:abscheme}
For any $0\leq a \leq b \leq T$ there is a randomised adaptive strategy for a cop such that the cop intercepts a random $R$-robber with probability at least $2^{-W(a,b)}$.
\end{lemma}
\begin{proof}
Let $(Y_s \colon s \in [s_0,s_1])$ be the predictions associated to the $(a,b)$-scheme $(s_0,s_1,\tau)$. Note that, by \ref{i:cstpred}, since $\tau_s=b$ is constant for all $ s \in [s_0,a]$, it follows that $Y_s = Y_{s_0}$ for all $s \in [s_0,a]$. Let us first show that we can define the strategy so that $C_a = Y_a= Y_{s_0}$. We split into two cases.

If $a < D$, and so $s_0 = -1$, then by \ref{i:initialY} the value $Y_{-1}$, which is sampled from the unconditional law of $R_b$, is revealed to the cop before he makes his initial placement. Hence, the cop can choose to start at the vertex $Y_{-1}$, and remain there until time $a$. Clearly, in this case, $C_a = C_0 = Y_{-1} = Y_a$.

If $a \geq D$, and so $s_0 = a-D$, the cop chooses his initial placement arbitrarily. At the start of round $a-D+1$ the cop has revealed the value of $Y_{s_0}$. The cop can then choose a walk of length $D$ from his current location to $Y_{s_0}$ and walk along it to arrive at $Y_{s_0}$ in round $a$. Hence, again in this case, we have arranged that $C_a = Y_{s_0} = Y_a$.

In each subsequent round $s+1$ with $a\leq s < s_1$, the cop reveals $Y_{s+1}$ and moves to its location, which is a legal cop move by \ref{i:feasible}, since $Y_{s+1}$ is drawn from $K_s(Y_s,\cdot)$. Any moves in subsequent rounds are made arbitrarily.

This strategy ensures that $C_s = Y_s$ for each $s \in [a,s_1]$. Note that, since $\tau_{s_1} \in \{s_1,s_1-1\}$, if $C_{s_1} = R_{\tau_{s_1}}$, then the cop catches the robber in round $s_1$ (either by moving onto the robber's location if $\tau_{s_1} = s_1-1$ or by the robber moving onto the cop's location if $\tau_{s_1} = s_1$). Hence, by \Cref{l:capture prob}, the cop catches the robber with probability at least
\[
\mathbb{P}(Y_{s_1} = R_{\tau_{s_1}})\geq 2^{-W(\tau)} = 2^{-W(a,b)}.
\]
\end{proof}

\begin{lemma}\label{l:capture}
For every random variable $R$ distributed on $\cc{W}_T$, there is some deterministic adaptive strategy for a cop who catches a random $R$-robber with probability at least $n^{-\frac{1}{2} - \frac{D}{2T+1}}$.
\end{lemma}
\begin{proof}
Let us consider the family of $(a,b)$-schemes where $(a,b) \in \cc{F}$ with
\[
\cc{F} = \{ (0,b) \colon 0 \leq b < T \} \cup \{ (a,T) \colon 0 \leq a\leq T\}
\]
where we note that $|\cc{F}| = 2T+1$. Let $W(\cc{F}) = \sum_{(a,b) \in \cc{F}} W(a,b)$. We will consider how many times each entry $\gamma_{s,t}$ of the information matrix appears in the sum corresponding to $W(\cc{F})$.

Looking at \eqref{e:locationandtrackingcost}, we see that the total cost, placement and interception, of the $(0,b)$-scheme consists of the sum of the entries of an anti-diagonal of the information matrix, specifically of $\{ \gamma_{s,t} \colon s\leq t, s+t =b\}$. Note that these are disjoint for distinct $0\leq b<  T$.

Similarly, if we look at the interception cost of the $(a,T)$-scheme we see it also consists of the entries of another anti-diagonal, this time $\{ \gamma_{s,t} \colon s\leq t <T ,s+t =a+T\}$. These are again disjoint for distinct $0 \leq a \leq T$, and disjoint from the previous anti-diagonals. Furthermore, between them these anti-diagonals cover every entry $\gamma_{s,t}$ of the information matrix with $s\leq t$ apart from the last column, the entries $\gamma_{s,T}$ with $0 \leq s \leq T$.

It remains then to consider the placement costs of the $(a,T)$-schemes, which consists of the sum of the entries $\{ \gamma_{s,T} \colon \max \{a-D +1, 0\} \leq s \leq a\}$. Note that every $\gamma_{s,T}$ is counted in at most $D$ such schemes, specifically for $s \leq a \leq \min \{T, s+D-1\}$.

Since by \eqref{e:chainrule} the entries in each column of the information matrix on or above the diagonal sum to at most $\log_2 n$, we see that
\[
W(\cc{F}) \leq \sum_{t=0}^{T-1} \sum_{s=0}^t \gamma_{s,t} + D \sum_{s=0}^T \gamma_{s,T} \leq (T+ D) \log_2 n.
\]
Hence, there is some $(a,b)$-scheme in $\cc{F}$ whose weight is at most the average
\[
W(a,b) \leq \frac{W(\cc{F})}{2T+1} \leq  \frac{T+D}{2T+1} \log_2 n \leq \left(\frac{1}{2} + \frac{D}{2T+1}\right) \log_2 n,
\]
and hence by \Cref{l:abscheme} there is a randomised strategy which catches a random $R$-robber with probability at least $n^{-\frac{1}{2} - \frac{D}{2T+1}}$. 

Since the cop's random choices are independent of the random trajectory $R$, his success probability is the average over these random choices of the success probabilities of the deterministic adaptive strategy given by these choices. Hence, there must exist some deterministic adaptive strategy whose success probability is at least the average, which is at least $n^{-\frac{1}{2} - \frac{D}{2T+1}}$, as claimed.
\end{proof}

\section{Bounding the cop number}
\begin{proof}[Proof of \Cref{t:copdiambound}]
Let $T= D \log_2 n$ and let $p=n^{-\frac{1}{2} - \frac{D}{2T+1}} = \Omega\left( n^{-\frac{1}{2}}\right)$.

Let us recursively define a sequence of sets $\cc{W}_T(G)=\cc{W}_1 \supseteq \cc{W}_2 \supseteq \ldots$ and random variables $R^1,R^2, \ldots$, where $R^i$ is uniformly distributed on $\cc{W}_i$. Note that $|\cc{W}_T(G)|\leq n(\Delta + 1)^T$.

Given $R^i$, by \Cref{l:capture} there is a deterministic strategy for a cop $c_i$ which catches a random $R^i$-robber with probability at least $p$, and hence intercepts at least $p |\cc{W}_i|$ of the possible robber trajectories in $\cc{W}_i$. We let $\cc{W}_{i+1}$ be the remaining walks.

We continue until the first $k$ such that $|\cc{W}_{k+1}| =0$, noting that, since $|\cc{W}_{i+1}| \leq (1-p)^i |\cc{W}_T(G)| \leq e^{-pi} |\cc{W}_T(G)|$, we have
\begin{align*}
k &\leq \frac{ 1 + \ln |\cc{W}_T(G)|}{p}\\
&= O \left( \sqrt{n} \left( \ln n + T \ln (\Delta+1) \right)\right)\\
&= O\left( T \sqrt{n} \ln (\Delta + 1) \right) \\
&=O\left( D \sqrt{n} \ln n \ln (\Delta + 1) \right),
\end{align*}
where we used that $n \leq (\Delta+1)^D \leq (\Delta + 1)^T$.

We claim that the cops $c_1, c_2, \ldots, c_k$ catch any robber, and hence $c(G) \leq k$, from which the result follows. Indeed, suppose for a contradiction that the robber survived to round $T$ and write $\bm{r} = (r_0,r_1, \ldots , r_T)$ for its realised trajectory. Since $\bm{r} \in \cc{W}_T(G)$, by construction there is some $i\in [k]$ such that $\bm{r} \in \cc{W}_i \setminus \cc{W}_{i+1}$. Then, by construction the cop $c_i$ intercepts the trajectory $\bm{r}$ by time $T$. Indeed, the cop $c_i$ has the same initial position, and sees the same sequence of prefixes $(r_0,r_1,\ldots,r_s)$ as when its strategy is played against the fixed walk $\bm{r}$ and hence makes the same moves and will intercept the robber by time $T$, contradicting our assumption.
\end{proof}

\begin{proof}[Proof of \Cref{t:weakmeyniel}]
Let $d= \frac{n^{\frac{1}{4}}}{\sqrt{\ln n \ln (\Delta +1)}}$. By a classic result of Aigner and Fromme \cite{AF84}, given any shortest path $P$ in a graph $H$, a single cop has a strategy to \emph{guard} $P$. That is, the cop can ensure that, after a finite amount of time, if the robber ever enters a vertex of $P$, then he is caught in the next round. Hence, we may assume that after this time the robber is restricted to the vertices and edges of the graph $H \setminus V(P)$.

Let us choose a maximal sequence of paths $P_1, P_2, \ldots , P_{k_1}$ such that each $P_i$ is a shortest path of length greater than $d$ in $G \setminus \left( \bigcup_{j=1}^{i-1} V(P_j)\right)$. Since each $P_i$ has length at least $d+1$ and the $P_i$ are disjoint, we have $k_1 \leq \frac{n}{d+1}$.

By the above, for each $i \in [k_1]$ there is a strategy for a single cop to guard $P_i$ in $G \setminus \left( \bigcup_{j=1}^{i-1} V(P_j)\right)$ and hence there is a strategy for $k_1$ cops to guard $\bigcup_{j=1}^{k_1} V(P_j)$. Hence, after a finite amount of time the robber is restricted to the subgraph $H = G\setminus \left( \bigcup_{j=1}^{k_1} V(P_j)\right)$, and in particular to one of its components, $H_1, H_2, \ldots , H_m$. It follows that, $c(G) \leq \frac{n}{d+1} + \max_{i \in [m]} c(H_i)$. Indeed, once the $k_1$ cops have restricted the robber to some $H_i$ we can move a reserve team of $c(H_i)$ cops to their starting positions in $H_i$ under an optimal strategy, and from then on these cops follow that strategy. If the robber ever leaves $H_i$ he is caught by one of the initial $k_1$ cops, and if he never leaves $H_i$, then he is eventually caught by one of the $c(H_i)$ cops playing in $H_i$.

Note that, by maximality of the path family, the diameter of every $H_i$ is at most $d$ and clearly every $H_i$ has at most $n$ vertices and maximum degree at most $\Delta$. Hence, by \Cref{t:copdiambound}, 
\[
c(G) \leq \frac{n}{d+1} + O\left( d \sqrt{n} \ln n \ln (\Delta +1) \right) = O\left( n^{\frac{3}{4}} \sqrt{\ln n \ln (\Delta +1)} \right).
\]
\end{proof}
\subsection*{Acknowledgements}
This research was funded in whole or in part by the Austrian Science Fund (FWF) [10.55776/P36131]. For open access purposes, the author has applied a CC BY public copyright license to any author accepted manuscript version arising from this submission.
\bibliographystyle{abbrv}
\bibliography{cop}
\end{document}